# Graph Generated Union-closed Families of Sets

Emanuel Knill

Los Alamos National Laboratory

Los Alamos, NM 87545

knill@lanl.gov

1992, revised June 1993

**Abstract**

Let $G$ be a graph with vertices $V$ and edges $E$. Let $\mathcal{F} = \{\bigcup_{e \in \mathcal{E}} e \mid \mathcal{E} \subseteq E\}$ be the union-closed family of sets generated by $E$. Then $\mathcal{F}$ is the family of subsets of $V$ without isolated points. Theorem: There is an edge $e \in E$ such that $\left|\{U \in \mathcal{F} \mid U \supseteq e\}\right| \leq \frac{1}{2}|\mathcal{F}|$. This is equivalent to the following assertion: If $\mathcal{H}$ is a union-closed family generated by a family of sets of maximum degree two, then there is an $x$ such that $\left|\{U \in \mathcal{H} \mid x \in U\}\right| \geq \frac{1}{2}|\mathcal{H}|$. This is a special case of the union-closed sets conjecture. To put this result in perspective, a brief overview of research on the union-closed sets conjecture is given. A proof of a strong version of the theorem on graph-generated families of sets is presented. This proof depends on an analysis of the local properties of $\mathcal{F}$ and an application of Kleitman's lemma. Much of the proof applies to arbitrary union-closed families and can be used to obtain bounds on $\left|\{U \in \mathcal{F} \mid U \supseteq e\}\right|/|\mathcal{F}|$.

## 1 Preliminaries

Unless stated otherwise, all sets are assumed to be finite. Let $\mathbf{N}$ be the set of non-negative integers. Define $[n] = \{i \in \mathbf{N} \mid 1 \leq i \leq n\}$ and $[m,n] = \{i \in$ $\mathbf{N} \mid m \leq i \leq n\}$. Let $X$ be a set. A *family* of sets on $X$ is a subset $\mathcal{F}$ of the power set $2^X$ of $X$. $\mathcal{F}$ is *union-closed* (*intersection-closed*) iff for every $U, V \in \mathcal{F}$, $U \cup V \in \mathcal{F}$ ($U \cap V \in \mathcal{F}$). The family $\mathcal{F}$ is a poset, where the members of $\mathcal{F}$ are ordered by inclusion. For $Y \subseteq X$, let

$$\mathcal{F}_{\subseteq Y} = \{U \in \mathcal{F} \mid U \subseteq Y\}$$

be the family *induced by $\mathcal{F}$ in $Y$*, let

$$\mathcal{F}_{\supseteq Y} = \{U \in \mathcal{F} \mid U \supseteq Y\}$$

be the family *induced by $\mathcal{F}$ above $Y$*, let

$$\mathcal{F}_{\cap Y} = \{U \cap Y \mid U \in \mathcal{F}\}$$

be the *restriction of $\mathcal{F}$ to $Y$* and let

$$\mathcal{F}_{\setminus Y} = \{U \setminus Y \mid U \in \mathcal{F}\}$$

be the *restriction of $\mathcal{F}$ away from $Y$*. If $Y = \{x\}$ is a singleton, the set brackets are omitted. The *degree* of $Y$ in $\mathcal{F}$ is defined as $d_\mathcal{F}(Y) = |\mathcal{F}_{\supseteq Y}|$. For families of sets $\mathcal{F}$ and $\mathcal{G}$, let $\mathcal{F} \vee \mathcal{G} = \{U \cup V \mid U \in \mathcal{F}, V \in \mathcal{G}\}$ and $\mathcal{F} \wedge \mathcal{G} = \{U \cap V \mid U \in \mathcal{F}, V \in \mathcal{G}\}$. The *transpose* or *dual* of $\mathcal{F}$ is the collection of sets $\mathcal{F}^*$ on $\mathcal{F}$ defined by

$$\mathcal{F}^* = \langle \mathcal{F}_{\supseteq x} \mid x \in X \rangle,$$

where the $\mathcal{F}_{\supseteq x}$ are counted with multiplicities. A collection of sets is called *simple*, if each member occurs only once (i.e. if it is a family of sets). If the dual of $\mathcal{F}$ is not simple, then there are elements of $X$ which are not separated by any member of $\mathcal{F}$. Call $\mathcal{F}$ *primitive* on $X$ if the dual of $\mathcal{F}$ is simple and $\bigcup \mathcal{F} = X$.

Let $P$ be a poset. The *dual* $P^*$ of $P$ is $P$ with the reverse ordering. A map $f$ from $P$ to a poset $Q$ is *order-preserving* iff $x \leq y$ implies $f(x) \leq f(y)$.

The set of order-preserving maps from $P$ to $Q$ is denoted by $Q^P$. $Q^P$ is a poset with the pointwise ordering $\pi \leq \sigma$ iff for all $x$, $\pi(x) \leq \sigma(x)$. The element $x \in P$ covers $y \in P$ iff $x > y$ and $x \geq z \geq y$ implies $z = x$ or $z = y$. For $A \subseteq P$, let

$$(A]_P = \{x \in P \mid x \leq y \text{ for some } y \in A\}$$
$$[A)_P = \{x \in P \mid x \geq y \text{ for some } y \in A\}.$$

Thus $(A]_P$ is the *order ideal* and $[A)_P$ is the *order filter* generated by $A$. Subscripts are omitted if the poset is clear from context. To simplify notation, for $x \in P$, let $[x) = [\{x\})$ and $(x] = (\{x\}]$.

If $P$ is closed under least upper bounds, then $P$ is a *join-semilattice* with the least upper bound of $x$ and $y$ denoted by $x \vee y$ (the *join* of $x$ and $y$). Similarly, if $P$ is closed under greatest lower bounds, then $P$ is a *meet-semilattice* with the greatest lower bound of $x$ and $y$ denoted by $x \wedge y$ (the *meet* of $x$ and $y$). A *lattice* is both a join-semilattice and a meet-semilattice. Recall that any finite join-semilattice can be made into a lattice by adding a least element $\hat{0}$ if necessary; and similarly for meet-semilattices. The meet (join) of elements in a join-semilattice (meet-semilattice) is understood to be the meet (join) in this extension.

If $L$ is a lattice and $x \in L$ satisfies that $u \vee v = x$ implies $u = x$ or $v = x$, then $x$ is *join-irreducible*. Every element of $L$ can be represented as the join of the join-irreducible elements below it. The set of join-irreducible elements excluding the least element is denoted by $J(L)$. *Meet-irreducible* elements are defined dually. $M(L)$ denotes the set of meet-irreducible elements of $L$ excluding the largest element.

There are canonical correspondences between finite semilattices and primitive union- or intersection-closed families of sets. Clearly every union- or intersection-closed family considered as a poset is a semilattice. Let $L$ be a

meet-semilattice. Let

$$\mathcal{F}(L) = \{(x]_L \cap J(L) \mid x \in L\}.$$

Then $\mathcal{F}$ is a primitive intersection-closed family of sets. If $L$ is a join-semilattice, then the corresponding union-closed family of sets consists of the family of complements in $M(L)$ of members of $\mathcal{F}(L^*)$.

Let $\mathcal{F}$ be a primitive union-closed family of sets on $X = \bigcup \mathcal{F}$. The meet-irreducible elements of $\mathcal{F}$ are given by the sets $M_x = \bigcup \mathcal{F}_{X\setminus\{x\}}$. The union generators of $\mathcal{F}$ are given by $G(\mathcal{F}) = J(\mathcal{F} \cup \{\emptyset\})$. This is a family of sets with the property that no member is the union of any collection of the other members. If $G(\mathcal{F}) = \mathcal{G}$, define $\mathcal{F}(\mathcal{G}) = \mathcal{F}$.

A family of sets $\mathcal{G}$ is a *graph* iff for every $U \in \mathcal{G}$, $1 \leq |U| \leq 2$. The members of $\mathcal{G}$ are referred to as *edges*. The graph $\mathcal{G}$ is simple iff for every $U \in \mathcal{G}$, $|U| = 2$. A *graph generated* union-closed family of sets satisfies that $G(\mathcal{F})$ is a graph.

## 2 Union-closed families of sets

The following problem has been attributed to P.Frankl (Duffus [4], Stanley [14]).

**Problem 2.1 The union-closed sets conjecture.** *Let $\mathcal{F}$ be a union-closed family of sets with at least one non-empty set. Does there always exist an element $x$ such that $d_{\mathcal{F}}(x) \geq |\mathcal{F}|/2$?*

The family of sets consisting of only the empty set does not satisfy the assertion of this problem. Henceforth all families of sets and semilattices are assumed to have at least two members. It is straightforward to check that the union-closed conjecture holds for $\mathcal{F} = 2^{[n]}$ if $n \geq 1$.

Using duality and the correspondence between union-closed families of sets and join-semilattices one arrives at other questions equivalent to the union-closed sets conjecture. The first such question is obtained by complementing the sets of a union-closed family. This yields an intersection-closed family order-isomorphic to the poset dual.

**Problem 2.2** *Let $\mathcal{F}$ be an intersection-closed family of sets. Does there always exist an element $x \in \bigcup \mathcal{F}$ such that $d_{\mathcal{F}}(x) \leq |\mathcal{F}|/2$?*

The second equivalent question is obtained from the first by using the correspondence between intersection-closed families and meet-semilattices and by observing that if $\mathcal{F}$ is a primitive intersection-closed family, then the elements of $\bigcup \mathcal{F}$ are in one-to-one correspondence with $J(\mathcal{F})$.

**Problem 2.3** *Let $L$ be a meet-semilattice. Does there always exist a member $x \in J(L)$ such that $|[x)| \leq |L|/2$?*

There is a corresponding version for join-semilattices. The last equivalent question to be given here is obtained by considering a union-closed family of sets with $\emptyset$ as a meet-semilattice.

**Problem 2.4** *Let $\mathcal{F}$ be a union-closed family of sets with $\emptyset \in \mathcal{F}$. Is there always a member $U$ of $G(\mathcal{F})$ such that $|\mathcal{F}_{\supseteq U}| \leq |\mathcal{F}|/2$?*

The union-closed sets conjecture can be generalized as follows (Knill [7]). Let $P$ be a poset and let $p$ be the number of filters of $P$. A meet-semilattice $L$ has the *P-density property* iff there exists $x \in J(L)$ such that $|[x)^P|/|L^P| \leq 1/p$. The quantity on the left of this inequality is called the *P-density* of $x$ in $L$. The *density* of $x$ is the [1]-density of $x$. The *density property* is the [1]-density property.

**Problem 2.5 The $P$-density problem.** *What meet semilattices have the $P$-density property?*

For $P = [1]$, this reduces to Problem 2.3. I.e. $L$ has the density property iff it satisfies the assertion of Problem 2.3. If $L$ has the density property with witness $a \in J(L)$, then we say that $(L, a)$ has the density property. Again the trivial example, the one-element semilattice, does not have the $P$-density property.

Some results related to the union-closed sets conjecture are given next. Recall that any intersection-closed family is a meet-semilattice.

The next two results have been rediscovered a few times and are given in [12] [10].

**Theorem 2.6** *If the intersection-closed family $\mathcal{F}$ contains $(\bigcup \mathcal{F}) \setminus \{x, y\}$ for some $x, y \in \bigcup \mathcal{F}$, then $\mathcal{F}$ has the density property.*

Let $S(\mathcal{F}) = \sum_{U \in \mathcal{F}} |U|$.

**Theorem 2.7** *If the average size $S(\mathcal{F})/|\mathcal{F}|$ of the members of the intersection-closed family $\mathcal{F}$ is at most $\frac{1}{2}|\bigcup \mathcal{F}|$, then $\mathcal{F}$ has the density property.*

Theorems 2.6 and 2.7 are obtained by double counting. They can be generalized substantially; a long overdue analysis of the technique is given by Poonen [10]. These ideas can be used to establish lower bounds on counter-examples to the union-closed sets conjecture (currently these bounds are about $|\mathcal{F}| \geq 29$ and $|\bigcup \mathcal{F}| \geq 8$).

The motivation underlying the work presented in Section 3 is related to the ideas in [10]. One of the main results of [10] (Theorem 1) characterizes union-closed families $\mathcal{F}$ which satisfy that for every union-closed extension $\mathcal{G} \supseteq \mathcal{F}$ there exists an $x \in \bigcup \mathcal{F}$ with $d_{\mathcal{G}}(x) \geq |\mathcal{G}|/2$. Here we consider extensions $\mathcal{G}$ of $\mathcal{F}$ in the context of Problem 2.4. In terms of the union-closed sets conjecture, this

can be seen to involve modifying $\mathcal{F}$ by adding new generators and adding new elements to the generators of $\mathcal{F}$.

For a given union-closed family $\mathcal{F}$ with $\bigcup \mathcal{F} \subseteq X$, define

$$s(\mathcal{F}, X) = \frac{S(\mathcal{F})}{|\mathcal{F}||X|}.$$

As in Theorem 2.7, if $s(\mathcal{F}, \bigcup \mathcal{F}) \geq \frac{1}{2}$ then $\mathcal{F}$ satisfies the union-closed sets conjecture. Wójcik [16] studies the properties of $s_m = \min\{s(\mathcal{F}, X) \mid |X| = m, \bigcup \mathcal{F} = X\}$. He conjectures that $s_m = (1 + o(1))\frac{\log_2 m}{2m}$ as $m \to \infty$. In fact, this is one of very few so far unproven, apparently *proper* consequences of the union-closed sets conjecture. Wójcik also gives an exact version of his conjecture which follows from a slight extension of the union-closed sets conjecture.

Wójcik's conjecture can be strengthened as follows (unpublished): For each $n \geq 0$, let $n = \sum_{i=0}^{l(n)} b(n)_i 2^i$ be the binary expansion of $n$, where $n = 0$ and $l(n) = 0$ or $b(n)_{l(n)} = 1$. Let

$$U(n) = \{i \geq 0 \mid i = l(n) \text{ and } b(n)_i = 1 \text{ or } i < l(n) \text{ and } b(n)_i = 0\}.$$

The function $U$ is a bijection from $\mathbf{N}$ to finite subsets of $\mathbf{N}$. It induces an ordering on such sets which is closely related to the lexicographic orderings. If $k < l$ and $U(i) = U(k) \cup U(l)$, then $i \leq l$. It follows that the families $\mathcal{U}(n) = \{U(0), U(1), \ldots, U(n-1)\}$ are union-closed. Let

$$t_n = \min\{S(\mathcal{F}) \mid |\mathcal{F}| = n\}. \tag{1}$$

**Problem 2.8** *Is it true that $t_n = S(\mathcal{U}(n))$?*

This can also be shown to be a consequence of the union-closed sets conjecture. It seems possible that the families $\mathcal{U}(n)$ minimize $\max_x d_{\mathcal{F}}(x)$ for union-closed families $\mathcal{F}$ with $|\mathcal{F}| = n$. This would imply the union-closed sets conjecture. Note that $\mathcal{F}(2^n)$ is the family of all subsets of $[0, n-1]$. Problem 2.8

is reminiscent of the Kruskal-Katona theorem according to which the sizes of the shadow and the ideal generated by a family of $n$ $k$-sets are minimized by the family which consists of the first $n$ $k$-sets in the squashed ordering (see [1], ch. 7).

Many of the standard classes of lattices can be shown to satisfy $P$-density properties. A lower-semimodular coatom of the lattice $L$ is a coatom $a$ with the property that for every $w \in L$ with $w \not\leq a$, $w$ covers $a \wedge w$.

**Theorem 2.9** *The following classes of lattices have the indicated density property.*

- *Distributive lattices ($P$-density for all $P$).*

- *Modular lattices ($P$-density for all $P$).*

- *Geometric lattices ($[n]$-density for all $n$; this is stated for $n = 1$ in [4], see [7]).*

- *Lattices where every $(x]$ is complemented (density property [10]).*

- *Lattices with a lower-semimodular coatom ($P$-density for all $P$ [7]).*

- *Lattices of height $h$ with $h = |J(L)|$ ($P$-density for all $P$).*

- *Selfdual lattices (density property).*

Note that the third class includes the second which includes the first, and the fifth class includes the first two. The result for selfdual lattices follows from the observation that for every lattice $L$ either $L$ or $L^*$ has the density property. The other results depend on the existence of certain matching properties in a lattice. If $F$ is a filter of $P$, then let $T(L, F, a) = \{\pi \in L^P \mid \pi(x) \geq a \text{ iff } x \in F\}$. The members of $T(L, F, a)$ are called the order-preserving maps of *type* $(F, a)$. We say that $L$ has the *P-matching property* iff there is a join-irreducible $a$ such that for every filter $F \subseteq P$, there is a decreasing one-to-one map $\rho : \mathrm{T}(L, P, a) \to \mathrm{T}(L, F, a)$. (A map $\rho$ is *decreasing* iff $\rho(x) \leq x$.) $L$ has the *full P-matching property* iff there is a join-irreducible $a$ such that for all filters $F \subseteq G \subseteq P$, there is a decreasing one-to-one map $\rho : \mathrm{T}(L, G, a) \to \mathrm{T}(L, F, a)$. The $P$-matching properties imply the $P$-density property. There are non-trivial lattices which do not satisfy any $P$-matching property. One such example is given by the union-closed family generated by the edges of the pentagon.

Here are some of the results about preservation of the density and matching properties under lattice operations obtained in [7].

**Theorem 2.10** *If the meet-semilattice $L$ has the $P$-density property, then so does $L \times M$ for any meet-semilattice $M$.*

For posets $P$ and $Q$, $P + Q$ is the disjoint union of $P$ and $Q$ ordered by the union of the orders on $P$ and $Q$.

**Theorem 2.11** *Preservation of matching properties for semilattices.*

- *If $L$ has the (full) $P$-matching property, then so does $L \times M$.*

- *If $L$ has the (full) $P$-matching property for $a \in J(L)$ and $I$ is an ideal of $L$ with $a \in I$, then $I$ has the (full) matching property for $a$.*

- *If $L$ has the (full) $P$-matching property, then so does $L^Q$.*

- *$L$ has the full $P + Q$-matching property for $a \in J(L)$ iff $L$ has the full $P$- and the full $Q$-matching property for $a$.*

- *If $L$ has the $P$- and the $Q$-matching property for $a \in J(L)$, then $L$ has the $P + Q$-matching property for $a$. The reverse implication holds if $a$ is an atom.*

Theorem 2.10 can be used to show that if the bound of Theorem 2.13 can be improved to $1 - \frac{1}{p} + o(1)$, then the statement of the $P$-density problem is true. In also suffices to show that the $P$-density property is satisfied by atomic lattices.

Matching properties are also preserved by certain subdirect products which preserve *local* properties of lattices. This notion of locality can be made precise by using the representation of lattices by union-closed families of sets. Let $L$ be a union-closed family of sets. The *lattice neighborhood* $\mathcal{N}_L(U)$ of $U \subseteq \bigcup J(L)$ is the union-closed family generated by $\{\emptyset\} \cup \{A \in J(L) \mid A \cap U \neq \emptyset\}$. Here is one of the results that can be obtained [7].

**Theorem 2.12** *If $x \in J(L)$ and $\mathcal{N}_L(x)$ is a geometric lattice, then $L$ has the $[n]$-density property (with witness $x$).*

If the union-closed sets conjecture is true, then for every union-closed family $\mathcal{F}$ of sets there exists a cover $Y$ of the non-empty members of $\mathcal{F}$ with at most $\lceil \log_2(|\mathcal{F}| + 1) \rceil$ elements. Recall that $Y$ is a cover of $\mathcal{G}$ iff $Y$ intersects every member of $\mathcal{G}$. Assuming the union-closed sets conjecture, $Y$ is obtained by first chosing $y_1$ such that $y_1$ covers at least half of the members of $\mathcal{F}$. Since $\mathcal{F}_2 = \mathcal{F}_{\subseteq \cup \mathcal{F} \setminus \{y_1\}}$ is union-closed, one can chose $y_2$ such that $y_2$ covers at least half of the members of $\mathcal{F}_2$. Continuing in this fashion, the desired cover $Y$ is obtained. The fact that such a set $Y$ exists can be shown without making use of the union-closed sets conjecture. It suffices to observe that for every minimal cover $Y$ of the non-empty members of $\mathcal{F}$, $\mathcal{F}_{\cap Y} \cup \{\emptyset\}$ is Boolean. This implies that $|Y| \leq \log_2(|\mathcal{F}| + 1)$. A related result gives an upper bound on the minimum $P$-density of join-irreducibles in a semilattice [7].

**Theorem 2.13** *The minimum $P$-density $\rho$ of a join-irreducible in the meet-semilattice $L$ satisfies $\rho \leq 1 - \frac{1}{\log_p(|L|)}$ asymptotically as $|L|$ goes to infinity.*

An indication that the bound of Theorem 2.13 is weak is given by the fact that it is asymptotically best possible for the generalization of meet-semilattices to multi-meet-semilattices, where every member other than $\hat{1}$ is assigned a multiplicity.

Before introducing the main result of the paper, here is another interesting result [7].

**Theorem 2.14** *For every fixed lattice $L$, there is an $N$ such that for all $n \geq N$, $L$ has the $[n]$-density property.*

This result is proved by considering the Zeta polynomial after eliminating the obvious cases and lattices of height $|J(L)|$.

## 3  Graph generated union-closed families of sets

The main result proven in this paper concerns union-closed families generated by a graph and Problem 2.4. The most easily stated version of the result is as follows.

**Theorem 3.1** *If $\mathcal{F}$ is a graph-generated union-closed family of sets with $\emptyset \in \mathcal{F}$, then there is a set $U \in G(\mathcal{F})$ such that $|\mathcal{F}_{\supseteq U}|/|\mathcal{F}| \leq 1/2$.*

The proof of Theorem 3.1 given here actually shows a substantially stronger result of a local nature. Instead of assuming that $G(\mathcal{F})$ is a graph, it suffices that there is a $U \in G(\mathcal{F})$ such that an appropriate neighborhood of $U$ is graph-generated and $U$ has locally minimal degree. The proof involves estimating $|\mathcal{F}_{\supseteq U}|/|\mathcal{F}|$ from local properties of $\mathcal{F}$. These estimates are very general and can be applied to arbitrary union-closed families.

**Definition.** Let $\mathcal{G}$ be a graph and $U \in \mathcal{G}$. The *degree* in $\mathcal{G}$ of $U$ is

$$\mathrm{d}_{\mathcal{G}}(U) = \big|\{V \in \mathcal{G} \setminus \{U\} \mid V \cap U \neq \emptyset\}\big|.$$

**Definition.** Let $\mathcal{F}$ be a union-closed family of sets. The *closure* in $\mathcal{F}$ of the set $X$ is given by

$$\pi_{\mathcal{F}}(X) = \bigcup(\mathcal{F}_{\subseteq X}) = \bigcup\{V \in G(\mathcal{F}) \,|\, V \subseteq X\}.$$

The set of *isolated* elements of $X$ is

$$\overline{\pi}_{\mathcal{F}}(X) = X \setminus \pi_{\mathcal{F}}(X).$$

In terms of hypergraphs, $\overline{\pi}_{\mathcal{F}}(X)$ is the set of isolated points of the hypergraph $(\mathcal{F}_{\subseteq X}, X)$. For any subfamily $\mathcal{H}$ of $\mathcal{F}$, the *density* of $\mathcal{H}$ in $\mathcal{F}$ is the ratio $|\mathcal{H}|/|\mathcal{F}|$.

**Observation 3.2** *Let $\mathcal{F}$ be a union-closed family of sets. The following are equivalent:*

*(i) $U \in \mathcal{F}$,*

*(ii) $\pi_{\mathcal{F}}(U) = U$,*

*(iii) $\overline{\pi}_{\mathcal{F}}(U) = \emptyset$.*

Let $\mathcal{F}$ be a union-closed family of sets such that $\emptyset \in \mathcal{F}$ and let $U \subseteq \bigcup \mathcal{F}$. Recall that $\mathcal{N}_{\mathcal{F}}(U)$ (the lattice neighborhood in $\mathcal{F}$ of $U$) is the union-closed family generated by the empty set and the members $V$ of $G(\mathcal{F})$ with $V \cap U \neq \emptyset$. The third lattice neighborhood is given by $\mathcal{N}^3_{\mathcal{F}}(U) = \mathcal{N}_{\mathcal{F}}(\bigcup \mathcal{N}_{\mathcal{F}}(U))$. (The second one, which consists of the sets of $\mathcal{F}$ included in $\bigcup \mathcal{N}_{\mathcal{F}(U)}$ will not be used here.)

We now develop the technique for estimating the density of $\mathcal{F}_{\supseteq U}$ in $\mathcal{F}$. This estimate depends only on the third lattice neighborhood $\mathcal{N}^3_{\mathcal{F}}(U)$.

**Definition.** Let $\mathcal{H}$ be a family of sets. The *density* in $\mathcal{H}$ of the set $X$ is

$$\rho_{\mathcal{H}}(X) = \frac{|\mathcal{H}_{\supseteq X}|}{|\mathcal{H}|}.$$

Thus $\rho_{\mathcal{H}}(X)$ is the density of $\mathcal{H}_{\supseteq X}$ in $\mathcal{H}$. The reciprocal of $\rho_{\mathcal{H}}(X)$ is denoted by $(\frac{1}{\rho})_{\mathcal{H}}(X)$. Let $\rho = \rho_{\mathcal{F}}(U)$.

**Definition.** The union-closed family $\mathcal{F}'$ is a *(conservative) extension* of $(\mathcal{F}, U)$ iff there is a non-empty union-closed family of sets $\mathcal{H}$ such that $(\bigcup \mathcal{H}) \cap U = \emptyset$ and $\mathcal{F}' = \mathcal{F} \vee \mathcal{H}$.

Associativity of $\vee$ for families of sets yields:

**Observation 3.3** *The extension relation is transitive; i.e. if $\mathcal{F}_1$ is an extension of $(\mathcal{F}, U)$ and $\mathcal{F}_2$ is an extension of $(\mathcal{F}_1, U)$, then $\mathcal{F}_2$ is an extension of $(\mathcal{F}, U)$.*

Since $\mathcal{F} = \mathcal{F} \vee \{\emptyset\}$:

**Observation 3.4** $\mathcal{F}$ *is an extension of* $(\mathcal{F}, U)$.

**Definition.** Let

$$\tfrac{1}{\overline{\rho}} = \inf\{(\tfrac{1}{\rho})_{\mathcal{F}'}(U) \mid \mathcal{F}' \text{ is an extension of } (\mathcal{F}, U)\}.$$

**Observation 3.5** $\tfrac{1}{\overline{\rho}} \leq \tfrac{1}{\rho}$.

**Observation 3.6** *If $U \in J(\mathcal{F})$ and $\tfrac{1}{\overline{\rho}} \geq 2$, then $(\mathcal{F}, U)$ has the density property.*

The goal is to find lower bounds on $\tfrac{1}{\overline{\rho}}$ in terms of the local properties of $\mathcal{F}$ at $U$. To this end, let $\mathcal{F}'$ be an (arbitrary) extension of $(\mathcal{F}, U)$ and consider

$$(\tfrac{1}{\rho})_{\mathcal{F}'}(U) = \frac{|\mathcal{F}'|}{|\mathcal{F}'_{\supseteq U}|}.$$

**Definition.** For $X \cap U = \emptyset$, let

$$T_{\mathcal{F}'}(X) = \{Y \subseteq U \mid X \cup Y \in \mathcal{F}'\}.$$

We have

$$|\mathcal{F}'| = \sum_{X \in \mathcal{F}'_{\setminus U}} |T_{\mathcal{F}'}(X)|$$

and
$$\left|\mathcal{F}'_{\supseteq U}\right| \leq \left|\mathcal{F}'_{\backslash U}\right|,$$
where equality holds if $U \in \mathcal{F}$ (since $\mathcal{F}'$ is closed under union with members of $\mathcal{F}$). This gives
$$(\tfrac{1}{\rho})_{\mathcal{F}'}(U) \geq \frac{\sum_{X \in \mathcal{F}'_{\backslash U}} \left|T_{\mathcal{F}'}(X)\right|}{\left|\mathcal{F}'_{\backslash U}\right|}. \qquad (*)$$

Note that if $U \in \mathcal{F}$, then this is an identity, so that $\frac{1}{\rho} \leq 2^{|U|}$. We will determine a lower bound for $\left|T_{\mathcal{F}'}(X)\right|$ which is independent of $\mathcal{F}'$.

**Lemma 3.7** *Let $X \in \mathcal{F}'_{\backslash U}$. If $Y \subseteq U$ satisfies*

(i) $\pi_{\mathcal{F}}(X \cup Y) \cap U = Y$,

(ii) $\pi_{\mathcal{F}}(X \cup Y) \supseteq \pi_{\mathcal{F}}(X \cup U) \setminus U$,

*then $X \cup Y \in \mathcal{F}'$.*

Conditions (i) and (ii) are independent of $\mathcal{F}'$. Observe that (i) is equivalent to
$$\overline{\pi}_{\mathcal{F}}(X \cup Y) \cap U = \emptyset$$
and (ii) is equivalent to
$$\overline{\pi}_{\mathcal{F}}(X \cup Y) \setminus U \subseteq \overline{\pi}_{\mathcal{F}}(X \cup U).$$

Lemma 3.10 below gives conditions equivalent to (i) and (ii) in terms of the family of generators of $\mathcal{F}$.

**Proof.** Since $X \in \mathcal{F}'_{\backslash U}$, there exists $Z \subseteq U$ such that $X \cup Z \in \mathcal{F}'$. Since $\mathcal{F}'$ is an extension of $(\mathcal{F}, U)$, there is a union-closed (non-empty) family $\mathcal{H}$ such that $(\bigcup \mathcal{H}) \cap U = \emptyset$ and $\mathcal{F}' = \mathcal{F} \vee \mathcal{H}$. We have $X \cup Z = A \cup B$ for some $A \in \mathcal{F}$ and $B \in \mathcal{H}$. We can assume that $A = \pi_{\mathcal{F}}(X \cup Z)$. Let $A' = \pi_{\mathcal{F}}(X \cup Y)$. By (i), $A' \cap U = Y$. If we can show that $X \setminus A' \subseteq B$, then $X \cup Y = A' \cup B \in \mathcal{F}'$ and we are done.

Let $x \in X \setminus A'$. Then $x \in \overline{\pi}_\mathcal{F}(X \cup Y) \setminus U$. By (ii), $x \in \overline{\pi}_\mathcal{F}(X \cup U) \setminus U$. Since $Z \subseteq U$, $\overline{\pi}_\mathcal{F}(X \cup U) \subseteq \overline{\pi}_\mathcal{F}(X \cup Z)$, which yields $x \in \overline{\pi}_\mathcal{F}(X \cup Z) \setminus U$. Since $\overline{\pi}_\mathcal{F}(X \cup Z) \setminus U = X \setminus A \subseteq B$, we have $x \in B$, as desired. ∎

**Definition.** For every set $X$ with $X \cap U = \emptyset$, let $E_{\mathcal{F},U}(X)$ consist of the subsets $Y$ of $U$ satisfying the conditions of Lemma 3.7:

$$E_{\mathcal{F},U}(X) = \{Y \subseteq U \mid \pi_\mathcal{F}(X \cup Y) \cap U = Y \text{ and }$$
$$\pi_\mathcal{F}(X \cup Y) \supseteq \pi_\mathcal{F}(X \cup U) \setminus U\}.$$

Let $E(X) = E_{\mathcal{F},U}(X)$.

**Example 3.8** Suppose that $\mathcal{F}$ is the union-closed family of sets generated by the edges of the graph depicted in Figure 1 and the empty set.

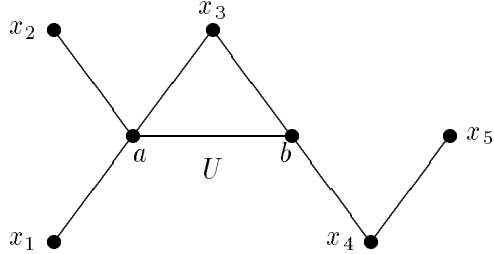

Figure 1

Then

$$J(\mathcal{F}) = \Big\{\{a,b\}, \{x_1,a\}, \{x_2,a\}, \{x_3,a\}, \{x_3,b\}, \{x_4,b\}, \{x_4,x_5\}\Big\}.$$

Let $U = \{a, b\}$. We have

$$\begin{aligned} E(\{x_4, x_5\}) &= \{\emptyset, \{b\}, \{a,b\}\}, \\ E(\{x_3\}) &= \{\{a\}, \{b\}, \{a,b\}\}, \\ E(\{x_1, x_4\}) &= \{\{a,b\}\}. \end{aligned}$$

By Lemma 3.7, $E(X) \subseteq T_{\mathcal{F}'}(X)$ for every $X \in \mathcal{F}'_{\setminus U}$. Using inequality $(*)$, we get

$$(\tfrac{1}{\rho})_{\mathcal{F}'}(U) \geq \frac{\sum_{X \in \mathcal{F}'_{\setminus U}} |E(X)|}{|\mathcal{F}'_{\setminus U}|}. \tag{$**$}$$

**Definition.** Let

$$\mu_{\mathcal{F}, \mathcal{F}'}(U) = \frac{\sum_{X \in \mathcal{F}'_{\setminus U}} |E(X)|}{|\mathcal{F}'_{\setminus U}|}.$$

Let $\mu_{\mathcal{F}'} = \mu_{\mathcal{F}, \mathcal{F}'}(U)$.

By inequality $(**)$:

**Observation 3.9** $\mu_{\mathcal{F}'} \leq (\tfrac{1}{\rho})_{\mathcal{F}'}(U)$.

We can define $E(X)$ in terms of the set of generators of $\mathcal{F}$.

**Lemma 3.10** *The set $Y$ is in $E(X)$ iff*

(i)′ *for every $x \in Y$ there exists $V \in J(\mathcal{F})$ with $x \in V \subseteq X \cup Y$,*

(ii)″ *for every $V \in J(\mathcal{F})$, if $V \setminus U \subseteq X$, then for every $x \in V \setminus U$ there exists $V' \in J(\mathcal{F})$ with $x \in V' \subseteq (X \cup Y)$.*

**Proof.** In fact, condition (i) of Lemma 3.7 is equivalent to (i)′ and condition (ii) of Lemma 3.7 is equivalent to (ii)″.

If $x \in W \in \mathcal{F}$, then there is a generator $V \in J(\mathcal{F})$ such that $x \in V \subseteq W$. This implies that $\pi_{\mathcal{F}}(X \cup Y) \cap U = Y$ iff (i)′ holds, hence (i) iff (i)′.

Suppose that $X$ and $Y$ satisfy (ii). Let $V \in J(\mathcal{F})$ and $V \setminus U \subseteq X$. Then $V \subseteq \pi_{\mathcal{F}}(X \cup U)$. By (ii), $\pi_{\mathcal{F}}(X \cup Y) \supseteq \pi_{\mathcal{F}}(X \cup U) \setminus U$. Hence, if $x \in V \setminus U$, then $x \in \pi_{\mathcal{F}}(X \cup Y)$, which implies that there is a generator $V'$ of $\mathcal{F}$ such that $x \in V' \subseteq X \cup Y$. Thus (ii)″ holds.

Conversely, suppose that $X$ and $Y$ satisfy (ii)″. We show that $\pi_{\mathcal{F}}(X \cup Y) \supseteq \pi_{\mathcal{F}}(X \cup U) \setminus U$. Let $x \in \pi_{\mathcal{F}}(X \cup U) \setminus U$. Then there exists $V \in J(\mathcal{F})$ such that $x \in V \subseteq X \cup U$. We have $V \setminus U \subseteq X$, so by (ii)″, there exists $V' \in J(\mathcal{F})$ with $x \in V' \subseteq X \cup Y$. This implies that $x \in \pi_{\mathcal{F}}(X \cup Y)$, as desired. ∎

**Definition.** The *neighborhood* in $\mathcal{F}$ of a set $X$ is given by

$$N_{\mathcal{F}}(X) = X \cup \bigcup \{V \in J(\mathcal{F}) \mid V \cap X \neq \emptyset\}.$$

Let $N = N_{\mathcal{F}}(U)$ and $N^2 = N_{\mathcal{F}}(N_{\mathcal{F}}(U))$. Observe that if $U \in \mathcal{F}$, then $N = \bigcup \mathcal{N}_{\mathcal{F}}(U)$ and $N^2 = \bigcup \mathcal{N}_{\mathcal{F}}^3(U)$.

**Lemma 3.11** *The family $\mathcal{F}$ is an extension of $(\mathcal{N}_{\mathcal{F}}^3(U), U)$.*

**Proof.** Let $\mathcal{H}$ be the union-closed family of sets generated by the empty set and the generators $V$ of $\mathcal{F}$ with $V \notin \mathcal{N}_{\mathcal{F}}^3(U)$. Then $\mathcal{F} = \mathcal{N}_{\mathcal{F}}^3(U) \vee \mathcal{H}$ and $(\bigcup \mathcal{H}) \cap U = \emptyset$. (This expresses $\mathcal{F}$ as an internal subdirect product of $\mathcal{N}_{\mathcal{F}}^3(U)$ and $\mathcal{H}$.) ∎

**Theorem 3.12** *Let $\mathcal{F}'$ be an extension of $(\mathcal{F}, U)$. Then*

$$\mu_{\mathcal{F},\mathcal{F}'}(U) = \mu_{\mathcal{N}_{\mathcal{F}}^3(U),\mathcal{F}'}(U).$$

**Proof.** By Lemma 3.11 and Observation 3.3, $\mathcal{F}'$ is an extension of $(\mathcal{N}_{\mathcal{F}}^3(U), U)$, so that $\mu_{\mathcal{N}_{\mathcal{F}}^3(U),\mathcal{F}'}(U)$ is well-defined. By Lemma 3.10 and by definition of $\mathcal{N}_{\mathcal{F}}^3(U)$, whether a given subset of $U$ is in $E(X)$ depends only on the generators $V \in J(\mathcal{F})$ with $V \in \mathcal{N}_\mathcal{F}^3(U)$. Since $J(\mathcal{N}_\mathcal{F}^3(U)) = J(\mathcal{F}) \cap \mathcal{N}_\mathcal{F}^3(U)$, Lemma 3.10 implies that

$$E_{\mathcal{F},U}(X) = E_{\mathcal{N}_\mathcal{F}^3(U),U}(X)$$

for every $X$ disjoint from $U$. The result now follows by definition of $\mu(U)$. ∎

**Definition.** Let $D \supseteq \bigcup \mathcal{F}$. An extension $\mathcal{F}'$ of $(\mathcal{F}, U)$ *minimizes* $\mu$ *in* $D$ iff $\bigcup \mathcal{F}' \subseteq D$ and for every extension $\mathcal{F}''$ of $(\mathcal{F}, U)$ with $\bigcup \mathcal{F}'' \subseteq D$, $\mu_{\mathcal{F}'} \leq \mu_{\mathcal{F}''}$.

**Theorem 3.13** *Let $D \supseteq \bigcup \mathcal{F}$. There exists an extension $\mathcal{F}' = \mathcal{F} \vee \mathcal{H}$ of $(\mathcal{F}, U)$ such that*

(i) *$\mathcal{H}$ is a filter of $2^{D \setminus U}$,*

(ii) *$\mathcal{F}'_{\setminus U} = \mathcal{H}$,*

(iii) *$\mathcal{F}'$ minimizes $\mu$ in $D$.*

**Proof.** Note that (i) implies (ii): Let $\mathcal{H}$ be a filter of $2^{D \setminus U}$. Since $\emptyset \in \mathcal{F}$, $(\mathcal{F} \vee \mathcal{H})_{\setminus U} \supseteq \mathcal{H}$. For the reverse inclusion, let $V \in \mathcal{F} \vee \mathcal{H}$. Then there is a $W \in \mathcal{H}$ with $W \subseteq V \setminus U$, hence $V \setminus U \in \mathcal{H}$.

Let $\mathcal{F} \vee \mathcal{G}$ be an extension of $(\mathcal{F}, U)$ such that $\bigcup \mathcal{G} \subseteq D \setminus U$ and $\mathcal{F} \vee \mathcal{G}$ minimizes $\mu$ in $D$. Since $(\mathcal{F} \vee (\mathcal{F} \vee \mathcal{G})_{\setminus U})_{\setminus U} = (\mathcal{F} \vee \mathcal{G})_{\setminus U}$, we have $\mu_{\mathcal{F} \vee (\mathcal{F} \vee \mathcal{G})_{\setminus U}} = \mu_{\mathcal{F} \vee \mathcal{G}}$. This implies that we can assume $\mathcal{G} = (\mathcal{F} \vee \mathcal{G})_{\setminus U}$.

Let $\mathcal{H}$ be the filter of $2^{D \setminus U}$ generated by $\mathcal{G}$ and let $\mathcal{F}' = \mathcal{F} \vee \mathcal{H}$. We show that $\mu_{\mathcal{F}'} \leq \mu_{\mathcal{F} \vee \mathcal{G}}$ which implies that $\mathcal{F}'$ is as desired.

For $X \in \mathcal{G}$, let

$$P(X) = \{Y \in \mathcal{H} \mid \pi_\mathcal{G}(Y) = X\}.$$

**Lemma 3.14** *Let $X, Y \in \mathcal{G}$ with $X \subseteq Y$. Then $|P(X)| \geq |P(Y)|$.*

**Proof.** Define for $Z \in P(Y)$

$$\sigma(Z) = (Z \setminus Y) \cup X.$$

Since every $Z \in P(Y)$ includes $Y$, $\sigma$ is a one-to-one map. To show that $\sigma$ maps $P(Y)$ into $P(X)$, let $Z \in P(Y)$. Let $X' = \pi_{\mathcal{G}}(\sigma(Z))$ and suppose that $X' \neq X$. Then $X' \supset X$. By definition of $\sigma(Z)$, $X' \not\subseteq Y$. Let $Y' = Y \cup X'$. Then $Y' \in \mathcal{G}$ and $Z \supseteq Y' \supset Y$, contradicting $Z \in P(Y)$. ∎

**Lemma 3.15** *Let $X \in \mathcal{G}$. If $Y \in P(X)$, then $E(Y) = E(X)$.*

**Proof.** We show that for every $Z \subseteq U$, $\pi_{\mathcal{F}}(X \cup Z) = \pi_{\mathcal{F}}(Y \cup Z)$. Let $Z \subseteq U$. Let $X' = \pi_{\mathcal{F}}(X \cup Z)$ and $Y' = \pi_{\mathcal{F}}(Y \cup Z)$. The inclusion $X \subseteq Y$ implies $X' \subseteq Y'$. Since $Y' \in \mathcal{F}$ and $X \in \mathcal{G}$, we have $Y' \cup X \in \mathcal{F} \vee \mathcal{G}$, so that $(Y' \cup X) \setminus U = (Y' \setminus U) \cup X \in \mathcal{G}$. Since $Y \in P(X)$ and $Y \supseteq (Y' \setminus U) \cup X$, it follows that $(Y' \setminus U) \cup X \subseteq \pi_{\mathcal{G}}(Y) = X$. Using $Y' \cap U \subseteq Z$, we get $Y' \subseteq X \cup Z$, hence $Y' \subseteq X'$.

The identities $\pi_{\mathcal{F}}(X \cup Z) = \pi_{\mathcal{F}}(Y \cup Z)$ and $\pi_{\mathcal{F}}(X \cup U) = \pi_{\mathcal{F}}(Y \cup U)$ imply that $Z \in E(X)$ iff $Z \in E(Y)$, as required. ∎

For $n \geq 0$, let

$$\mathcal{G}_{\leq n} = \{X \in \mathcal{G} \mid |P(X)| \leq n\}.$$

By Lemma 3.14, the $\mathcal{G}_{\leq n}$ are filters of $\mathcal{G}$. Since $(\mathcal{F} \vee \mathcal{G})_{\setminus U} = \mathcal{G}$, this implies that if $\mathcal{G}_{\leq n} \neq \emptyset$, then $\mathcal{F} \vee \mathcal{G}_{\leq n}$ is an extension of $(\mathcal{F}, U)$ such that $(\mathcal{F} \vee \mathcal{G}_{\leq n})_{\setminus U} = \mathcal{G}_{\leq n}$. Let $N$ be the maximum value of $|P(X)|$. Then $\mathcal{G}_{\leq N} = \mathcal{G}$. We use Lemmas 3.14 and 3.15 and the fact that the family $\{P(X) \mid X \in \mathcal{G}\}$ is a partition of $\mathcal{H}$ to compute $\mu_{\mathcal{F}'}$:

$$\mu_{\mathcal{F}'} = \frac{\sum_{X \in \mathcal{F}'_{\setminus U}} |E(X)|}{|\mathcal{F}'_{\setminus U}|}$$

$$\begin{aligned}
&= \frac{\sum_{X \in \mathcal{H}} |E(X)|}{\sum_{X \in \mathcal{H}} 1} \\
&= \frac{\sum_{X \in \mathcal{G}} |P(X)| \cdot |E(X)|}{\sum_{X \in \mathcal{G}} |P(X)|} \\
&= \frac{\sum_{n=0}^{N-1} \sum_{X \in \mathcal{G} \setminus \mathcal{G}_{\leq n}} |E(X)|}{\sum_{n=0}^{N-1} \sum_{X \in \mathcal{G} \setminus \mathcal{G}_{\leq n}} 1} \\
&= \frac{\sum_{n=0}^{N-1} \left( \sum_{X \in \mathcal{G}} |E(X)| - \sum_{X \in \mathcal{G}_{\leq n}} |E(X)| \right)}{\sum_{n=0}^{N-1} |\mathcal{G}| - |\mathcal{G}_{\leq n}|}.
\end{aligned}$$

Define $\mu_{\mathcal{F} \vee \emptyset} = \mu_{\mathcal{F} \vee \mathcal{G}}$. We have $\mu_{\mathcal{F} \vee \mathcal{G}} = \left( \sum_{X \in \mathcal{G}} |E(X)| \right) / |\mathcal{G}|$ and if $\mathcal{G}_{\leq n} \neq \emptyset$, then $\mu_{\mathcal{F} \vee \mathcal{G}_{\leq n}} = \left( \sum_{X \in \mathcal{G}_{\leq n}} |E(X)| \right) / |\mathcal{G}_{\leq n}|$ and $\mu_{\mathcal{F} \vee \mathcal{G}} \leq \mu_{\mathcal{F} \vee \mathcal{G}_{\leq n}}$. This yields

$$\begin{aligned}
\mu_{\mathcal{F}'} &= \frac{\sum_{n=0}^{N-1} \left( \mu_{\mathcal{F} \vee \mathcal{G}} |\mathcal{G}| - \mu_{\mathcal{F} \vee \mathcal{G}_{\leq n}} |\mathcal{G}_{\leq n}| \right)}{\sum_{n=0}^{N-1} \left( |\mathcal{G}| - |\mathcal{G}_{\leq n}| \right)} \\
&\leq \frac{\sum_{n=0}^{N-1} \left( \mu_{\mathcal{F} \vee \mathcal{G}} |\mathcal{G}| - \mu_{\mathcal{F} \vee \mathcal{G}} |\mathcal{G}_{\leq n}| \right)}{\sum_{n=0}^{N-1} \left( |\mathcal{G}| - |\mathcal{G}_{\leq n}| \right)} \\
&= \mu_{\mathcal{F} \vee \mathcal{G}} \frac{\sum_{n=0}^{N-1} \left( |\mathcal{G}| - |\mathcal{G}_{\leq n}| \right)}{\sum_{n=0}^{N-1} \left( |\mathcal{G}| - |\mathcal{G}_{\leq n}| \right)} \\
&= \mu_{\mathcal{F} \vee \mathcal{G}},
\end{aligned}$$

which implies that $\mathcal{F}'$ satisfies (iii), as required. ∎

**Lemma 3.16** *Let $X, Y$ be sets disjoint from $U$. If $Y \cap N^2 = X \cap N^2$, then $E(X) = E(Y)$.*

**Proof.** By Lemma 3.10, whether a given subset of $U$ is in $E(Z)$ depends only on the generators of $\mathcal{F}$ included in $N^2$. This implies that $E(Z)$ depends only on $Z \cap N^2$. ∎

We strengthen Theorem 3.13.

**Theorem 3.17** *Let $D \supseteq \bigcup \mathcal{F}$. There is an extension $\mathcal{F}' = \mathcal{F} \vee \mathcal{H}$ of $(\mathcal{F}, U)$ such that*

(i) $\mathcal{H}$ *is a filter of* $2^{D \setminus U}$,

(ii) $\mathcal{F}'_{\setminus U} = \mathcal{H}$,

(iii) $\mathcal{F}'$ *minimizes $\mu$ in $D$,*

(iv) $\mathcal{F}'_{\setminus N^2} = \{D \setminus N^2\}$.

**Proof.** By Theorem 3.13 there is an extension $\mathcal{F} \vee \mathcal{G}$ of $\mathcal{F}$ which satisfies (i), (ii) and (iii). Let $\mathcal{H} = \mathcal{G} \vee \{D \setminus N^2\}$. We show that $\mathcal{F}' = \mathcal{F} \vee \mathcal{H}$ is as required. By construction, $\mathcal{F}'$ satisfies (i), (ii) and (iv). For $X \in \mathcal{H}$, let

$$P(X) = \{Y \in \mathcal{G} \mid Y \cap N^2 = X \cap N^2\}.$$

Then $\{P(X) \mid X \in \mathcal{H}\}$ is a partition of $\mathcal{G}$. If $Y \in P(X)$, then by Lemma 3.16, $E(Y) = E(Y \cap N^2) = E(X \cap N^2) = E(X)$. Since $\mathcal{G}$ is a filter of $2^{D \setminus U}$, if $X \subseteq Y$, then $P(X)_{\setminus N^2} \subseteq P(Y)_{\setminus N^2}$. Since $P(X)_{\cap N^2} = \{X \cap N^2\}$ for each $X \in \mathcal{H}$, this implies that if $X \subseteq Y$, then $|P(X)| \leq |P(Y)|$. Let

$$\mathcal{H}_{\geq n} = \{X \in \mathcal{H} \mid |P(X)| \geq n\}.$$

Then $\mathcal{H}_{\geq n}$ is a filter included in $\mathcal{H}$ for each $n$ and $\mathcal{H}_{\geq 1} = \mathcal{H}$. Let $N$ be the maximum value of $|P(X)|$. By assumption on $\mathcal{F} \vee \mathcal{G}$, $\mu_{\mathcal{F} \vee \mathcal{H}_{\geq n}} \geq \mu_{\mathcal{F} \vee \mathcal{G}}$. Suppose that for some $n$, $\mu_{\mathcal{F} \vee \mathcal{H}_{\geq n}} > \mu_{\mathcal{F} \vee \mathcal{G}}$. Then

$$\begin{aligned}
\mu_{\mathcal{F} \vee \mathcal{G}} &= \frac{\sum_{X \in \mathcal{G}} |E(X)|}{|\mathcal{G}|} \\
&= \frac{\sum_{X \in \mathcal{H}} |P(X)| \cdot |E(X)|}{\sum_{X \in \mathcal{H}} |P(X)|} \\
&= \frac{\sum_{n=1}^{N} \sum_{X \in \mathcal{H}_{\geq n}} |E(X)|}{\sum_{n=1}^{N} \sum_{X \in \mathcal{H}_{\geq n}} 1}
\end{aligned}$$

$$
\begin{aligned}
&= \frac{\sum_{n=1}^{N}\left(|\mathcal{H}_{\geq n}|\mu_{\mathcal{F}\vee\mathcal{H}_{\geq n}}\right)}{\sum_{n=1}^{N}|\mathcal{H}_{\geq n}|} \\
&> \frac{\sum_{n=1}^{N}\left(|\mathcal{H}_{\geq n}|\mu_{\mathcal{F}\vee\mathcal{G}}\right)}{\sum_{n=1}^{N}|\mathcal{H}_{\geq n}|} \\
&= \mu_{\mathcal{F}\vee\mathcal{G}} \cdot \frac{\sum_{n=1}^{N}|\mathcal{H}_{\geq n}|}{\sum_{n=1}^{N}|\mathcal{H}_{\geq n}|} \\
&= \mu_{\mathcal{F}\vee\mathcal{G}},
\end{aligned}
$$

which is impossible. Hence, for each $n$, $\mathcal{F} \vee \mathcal{H}_{\geq n}$ also minimizes $\mu$ in $D$. In particular, $\mathcal{F}'$ minimizes $\mu$ in $D$, as required. ∎

**Corollary 3.18** *There is an extension $\mathcal{F}' = \mathcal{N}_{\mathcal{F}}^3(U) \vee \mathcal{H}$ of $(\mathcal{N}_{\mathcal{F}}^3(U), U)$ such that*

(i) *$\mathcal{H}$ is a filter of $2^{N^2 \setminus U}$,*

(ii) *$\mathcal{F}'_{\setminus U} = \mathcal{H}$,*

(iii) *if $\mathcal{F}''$ is an extension of $(\mathcal{F}, U)$, then $\mu_{\mathcal{F}, \mathcal{F}''}(U) \geq \mu_{\mathcal{N}_{\mathcal{F}}^3(U), \mathcal{F}'}(U)$.*

**Proof.** Let $D \supseteq \bigcup \mathcal{F}$ and let $\mathcal{F} \vee \mathcal{G}$ be an extension of $(\mathcal{F}, U)$ satisfying the conditions given in Theorem 3.17. By Theorem 3.12, $\mathcal{F} \vee \mathcal{G}$ is an extension of $(\mathcal{N}_{\mathcal{F}}^3(U), U)$ and

$$\mu_{\mathcal{F}, \mathcal{F} \vee \mathcal{G}}(U) = \mu_{\mathcal{N}_{\mathcal{F}}^3(U), \mathcal{F} \vee \mathcal{G}}(U).$$

By condition (iv) of Theorem 3.17, $(\mathcal{F} \vee \mathcal{G})_{\cap N^2}$ is isomorphic to $\mathcal{F} \vee \mathcal{G}$. This and Lemma 3.16 imply that

$$\mu_{\mathcal{N}_{\mathcal{F}}^3(U), \mathcal{F} \vee \mathcal{G}}(U) = \mu_{\mathcal{N}_{\mathcal{F}}^3(U), (\mathcal{F} \vee \mathcal{G})_{\cap N^2}}(U).$$

Since $\mathcal{F}$ is an extension of $(\mathcal{N}_{\mathcal{F}}^3(U), U)$, $\mathcal{F} = \mathcal{N}_{\mathcal{F}}^3(U) \vee \mathcal{G}'$ for some (non-empty) union-closed family $\mathcal{G}'$ with $\bigcup \mathcal{G}' \cap U = \emptyset$. The family $\mathcal{H}' = \mathcal{G}' \vee \mathcal{G}$ is a filter of $2^{D\setminus U}$ with $\mathcal{H}'_{\setminus N^2} = \{D \setminus N^2\}$. Let $\mathcal{H} = \mathcal{H}'_{\cap N^2}$. Then

$$\begin{aligned}(\mathcal{F} \vee \mathcal{G})_{\cap N^2} &= (\mathcal{N}^3_{\mathcal{F}}(U) \vee \mathcal{G}' \vee \mathcal{G})_{\cap N^2} \\ &= \mathcal{N}^3_{\mathcal{F}}(U) \vee \mathcal{H},\end{aligned}$$

so that $\mathcal{F}' = \mathcal{N}^3_{\mathcal{F}}(U) \vee \mathcal{H}$ is the required extension. ∎

Corollary 3.18 implies that to determine the minimum possible value of $\mu_{\mathcal{F},\mathcal{F}'}(U)$, $\mathcal{F}$ can be replaced by $\mathcal{N}^3_{\mathcal{F}}(U)$; thus assume that $\mathcal{F} = \mathcal{N}^3_{\mathcal{F}}(U)$.

Henceforth we assume that $U \in J(\mathcal{F})$. This implies that for every $X$ disjoint from $U$, $U \in E(X)$, so that $|E(X)| \geq 1$. To show that $(\mathcal{F}, U)$ has the density property, it suffices to show that for every extension $\mathcal{F}'$ of $(\mathcal{F}, U)$, $\mu_{\mathcal{F}'} \geq 2$.

**Theorem 3.19** *The following are equivalent:*

(i) *There exists an extension $\mathcal{F}'$ of $(\mathcal{F}, U)$ such that $\mu_{\mathcal{F}'} < 2$.*

(ii) *There is a filter $\mathcal{H}$ of $2^{N^2\setminus U}$ such that $\mu_{\mathcal{F} \vee \mathcal{H}} < 2$ and for every minimal member $X$ of $\mathcal{H}$, $|E(X)| = 1$.*

**Proof.** Assertion (ii) implies (i). Suppose that (i) holds. By Corollary 3.18, there is an extension $\mathcal{F} \vee \mathcal{H}$ of $(\mathcal{F}, U)$ such that $\mathcal{H}$ is a filter of $2^{N^2\setminus U}$ and $\mu_{\mathcal{F} \vee \mathcal{H}} < 2$. Let $\mathcal{H}$ be a minimal filter of $2^{N^2\setminus U}$ such that $\mu_{\mathcal{F} \vee \mathcal{H}} < 2$. Suppose that there is a minimal member $X$ of $\mathcal{H}$ such that $E(X) \geq 2$. Let $\mathcal{H}' = \mathcal{H}\setminus\{X\}$. Note that the assumption on $X$ and $\mu_{\mathcal{F} \vee \mathcal{H}} < 2$ imply that $\mathcal{H}' \neq \emptyset$. The family $\mathcal{H}'$ is a filter of $2^{N^2\setminus U}$ and

$$\begin{aligned}\mu_{\mathcal{F} \vee \mathcal{H}'} &= \frac{\sum_{Y \in \mathcal{H}'}|E(Y)|}{|\mathcal{H}'|} \\ &= \frac{\left(\sum_{Y \in \mathcal{H}}|E(Y)|\right) - |E(X)|}{|\mathcal{H}| - 1} \\ &< \frac{\sum_{Y \in \mathcal{H}}|E(Y)|}{|\mathcal{H}|}\end{aligned}$$

$$= \mu_{\mathcal{F} \vee \mathcal{H}}$$

$$< 2,$$

where we used the fact that if $a, b > 0$, $c > 1$ and $a/c < b$, then $(a-b)/(c-1) < a/c$. This contradicts the minimality assumption on $\mathcal{H}$, so that $\mathcal{H}$ is as desired. ∎

If there is no filter $\mathcal{H}$ satisfying the conditions in assertion (ii) of Theorem 3.19, then $\mu_{\mathcal{F}'}(U) \geq 2$ for every extension $\mathcal{F}'$ of $(\mathcal{F}, U)$, so that $(\mathcal{F}, U)$ has the density property. This yields:

**Observation 3.20** *If every extension $\mathcal{F} \vee \mathcal{H}$ of $(\mathcal{F}, U)$ such that*

(i) *$\mathcal{H}$ is a filter of $2^{N^2 \setminus U}$,*

(ii) *for every minimal $X \in \mathcal{H}$, $|E(X)| = 1$*

*satisfies $\mu_{\mathcal{F} \vee \mathcal{H}} \geq 2$, then $(\mathcal{F}, U)$ has the density property.*

Let $\mathcal{H}$ be an arbitrary filter satisfying (i) and (ii) of Observation 3.20 and let $\mathcal{F}' = \mathcal{F} \vee \mathcal{H}$. So far the discussion did not require any additional assumptions on $J(\mathcal{F})$. We now assume that $J(\mathcal{F})$ is a graph. Thus, for some $a, b \in \bigcup \mathcal{F}$, $U = \{a, b\}$. Assume that $a \neq b$ and let

$$N_a = \{x \in N \setminus U \mid \{x, a\} \in J(\mathcal{F}) \text{ and } \{x, b\} \notin J(\mathcal{F})\},$$
$$N_b = \{x \in N \setminus U \mid \{x, a\} \notin J(\mathcal{F}) \text{ and } \{x, b\} \in J(\mathcal{F})\},$$
$$N_{ab} = \{x \in N \setminus U \mid \{x, a\} \in J(\mathcal{F}) \text{ and } \{x, b\} \in J(\mathcal{F})\}.$$

Then $N \setminus U = N_a \cup N_b \cup N_{ab}$ (see Figure 2).

**Lemma 3.21** *If $X \in \mathcal{H}$, then $X \cap N_a \neq \emptyset$ and $X \cap N_b \neq \emptyset$.*

**Proof.** Let $X \in \mathcal{H}$. Let $X'$ be a minimal member of $\mathcal{H}$ such that $X' \subseteq X$. Suppose $X' \cap N_a = \emptyset$. Then either $\{b\} \in E(X')$ (if $X' \cap (N_{ab} \cup N_b) \neq \emptyset$) or $\emptyset \in E(X')$ (if $X' \cap (N_{ab} \cup N_b) = \emptyset$). Since $\{a,b\} \in E(X')$, this contradicts $|E(X')| = 1$. Thus $X' \cap N_a \neq \emptyset$. By symmetry, $X' \cap N_b \neq \emptyset$ and we are done. ∎

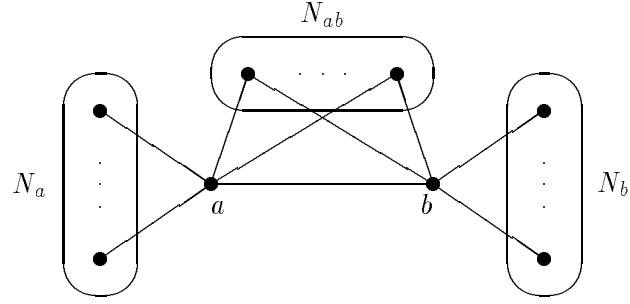

Figure 2

**Corollary 3.22** *If $X \in \mathcal{H}$ and $Y \subseteq U$, then $\pi_{\mathcal{F}}(X \cup Y) \cap U = Y$.*

**Proof.** Let $X \in \mathcal{H}$ and $Y \subseteq U$. Suppose that $a \in Y$. By Lemma 3.21, there exists $x \in N_a \cap X$, so that $\{x,a\} \in \mathcal{F}$ and $\{x,a\} \subseteq X \cup Y$, which implies that $a \in \pi_{\mathcal{F}}(X \cup Y)$. Similarly, if $b \in Y$, then $b \in \pi_{\mathcal{F}}(X \cup Y)$, as required. ∎

**Definition.** Let $Y \subseteq U$ and $x \in N \setminus U$. Let $\mathcal{E}(Y,x)$ consist of the subsets $X$ of $N^2 \setminus U$ such that there is an edge $\{x,y\} \in J(\mathcal{F})$ with $y \in X \cup Y$ or $y = x$. Define
$$\mathcal{E}(Y) = \bigcap_{x \in N \setminus U} \mathcal{E}(Y,x).$$

**Observation 3.23** *For $Y \subseteq U$ and $x \in N \setminus U$, $\mathcal{E}(Y,x)$ is a filter of $2^{N^2 \setminus U}$.*

**Observation 3.24** *If $X \in \mathcal{E}(Y)$, then $X$ and $Y$ satisfy condition (ii)″ of Lemma 3.10.*

**Lemma 3.25** *If $X \in \mathcal{E}(Y) \cap \mathcal{H}$, then $Y \in E(X)$.*

**Proof.** Let $X \in \mathcal{E}(Y) \cap \mathcal{H}$. By Observation 3.24, $X$ and $Y$ satisfy (ii)″ of Lemma 3.10. By Corollary 3.22, $X$ and $Y$ satisfy (i) of Lemma 3.7. The result follows by the proof of Lemma 3.10. ∎

**Definition.** Let $\mathcal{U}$ be a filter of $2^S$. Define
$$\nu(\mathcal{U}) = \frac{|\mathcal{U}|}{2^{|S|}}.$$
Thus $\nu(\mathcal{U})$ is the density of $\mathcal{U}$ in $2^S$.

**Theorem 3.26**
$$\mu_{\mathcal{F},\mathcal{F}'}(U) \geq 1 + \sum_{Y \subset U} \prod_{x \in N \setminus U} \nu\left(\mathcal{E}(Y,x)\right).$$

**Proof.** Let us start be recalling Kleitman's lemma. This result has many applications. Its proof can be found in Anderson [1].

**Theorem 3.27** (Kleitman [6]) *Let $X$ be an $n$-set. If $\mathcal{F}_1$ and $\mathcal{F}_2$ are filters in $2^X$, then*
$$\frac{|\mathcal{F}_1 \cap \mathcal{F}_2|}{2^n} \geq \frac{|\mathcal{F}_1|}{2^n} \frac{|\mathcal{F}_2|}{2^n}.$$

This says that the density (in $2^X$) of the intersection of two filters is at least the product of the densities of each.

Returning to the proof of Theorem 3.26, it follows from Kleitman's lemma that for every $Y \subseteq U$,
$$\nu(\mathcal{H} \cap \mathcal{E}(Y)) \geq \nu(\mathcal{H})\nu(\mathcal{E}(Y)).$$
Lemma 3.25 and the fact that $U \in E(X)$ for every $X \in \mathcal{H}$ yield
$$\begin{aligned}
\sum_{X \in \mathcal{H}} |E(X)| &= \sum_{Y \subseteq U} |\{X \in \mathcal{H} \mid Y \in E(X)\}| \\
&= |\{X \in \mathcal{H} \mid U \in E(X)\}| + \sum_{Y \subset U} |\{X \in \mathcal{H} \mid Y \in E(X)\}| \\
&\geq |\mathcal{H}| + \sum_{Y \subset U} |\mathcal{H} \cap \mathcal{E}(Y)|.
\end{aligned}$$

Using $\mathcal{H} = \mathcal{F}'_{\setminus U}$ we get

$$\begin{aligned}
\mu_{\mathcal{F}'} &= \frac{\sum_{X \in \mathcal{H}} |E(X)|}{|\mathcal{H}|} \\
&\geq \frac{|\mathcal{H}| + \sum_{Y \subset U} |\mathcal{H} \cap \mathcal{E}(Y)|}{|\mathcal{H}|} \\
&= 1 + \sum_{Y \subset U} \frac{\nu(\mathcal{H} \cap \mathcal{E}(Y))}{\nu(\mathcal{H})} \\
&\geq 1 + \sum_{Y \subset U} \nu(\mathcal{E}(Y)).
\end{aligned}$$

Multiple applications of Kleitman's Lemma yield

$$\nu(\mathcal{E}(Y)) \geq \prod_{x \in N \setminus U} \nu(\mathcal{E}(Y, x))$$

and the result follows. ∎

Define

$$J'(\mathcal{F}) = \{V \in J(\mathcal{F}) \mid |V| = 2\}.$$

For $x \in N \setminus U$, define

$$A(x) = \{\{x, y\} \mid y \notin U \text{ and } \{x, y\} \in J'(\mathcal{F})\}.$$

Let $n(x) = |A(x)|$.

**Lemma 3.28** *Let $x \in N \setminus U$. Then*

(i) $\nu(\mathcal{E}(\emptyset, x)) \geq 1 - (\frac{1}{2})^{n(x)}$.

(ii) *If $x \in N_a$, then $\nu(\mathcal{E}(\{a\}, x)) = 1$ and $\nu(\mathcal{E}(\{b\}, x)) \geq 1 - (\frac{1}{2})^{n(x)}$.*

(iii) *If $x \in N_b$, then $\nu(\mathcal{E}(\{a\}, x)) \geq 1 - (\frac{1}{2})^{n(x)}$, and $\nu(\mathcal{E}(\{b\}, x)) = 1$.*

(iv) *If $x \in N_{ab}$, then $\nu(\mathcal{E}(\{a\}, x)) = \nu(\mathcal{E}(\{b\}, x)) = 1$.*

**Proof.** Consider $\nu(\mathcal{E}(\emptyset, x))$. If $\{x\} \in J(\mathcal{F})$, then $\nu(\mathcal{E}(\emptyset, x)) = 1$. If $\{x\} \notin J(\mathcal{F})$, then

$$2^{N^2 \setminus U} \setminus \mathcal{E}(\emptyset, x) = 2^{N^2 \setminus (U \cup A(x))}.$$

We have
$$\left|2^{N^2\setminus(U\cup A(x))}\right| = 2^{|N^2\setminus U|}2^{-|A(x)|},$$
which implies that
$$\nu(\mathcal{E}(\emptyset,x)) \geq 1 - (\tfrac{1}{2})^{n(x)}.$$

The remaining cases are proved similarly. ∎

Recall that $d_{\mathcal{G}}(U)$ is the degree in the graph $\mathcal{G}$ of $U$. If $U$ has minimal degree in $J'(\mathcal{F})$, then $(\mathcal{F}, U)$ has the density property:

**Theorem 3.29** *If for every edge $V \in J'(\mathcal{F})$ with $V \cap U \neq \emptyset$, $d_{J'(\mathcal{F})}(V) \geq d_{J'(\mathcal{F})}(U)$, then for every extension $\mathcal{F}'$ of $(\mathcal{F}, U)$, $\mu_{\mathcal{F}, \mathcal{F}'}(U) \geq 2$.*

**Proof.** We can assume that $\mathcal{F}' = \mathcal{F} \vee \mathcal{H}$ where $\mathcal{H}$ is a filter satisfying (i) and (ii) of Observation 3.20. Let $n = d_{J'(\mathcal{F})}(U)$, $n_a = |N_a|$, $n_b = |N_b|$ and $n_{ab} = |N_{ab}|$. Then
$$n = n_a + n_b + 2n_{ab}.$$

Let $x \in N_a$. We have $d_{J'(\mathcal{F})}(\{a,x\}) \geq n$ and $d_{J'(\mathcal{F})}(\{a,x\}) = n_a + n_{ab} + n(x)$. This gives $n(x) \geq n - n_a - n_{ab} = n_b + n_{ab}$. Similarly, if $x \in N_b$, then $n(x) \geq n_a + n_{ab}$. Lemma 3.21 and the assumptions on $\mathcal{H}$ imply that $n_a \geq 1$ and $n_b \geq 1$. By Theorem 3.26 and Lemma 3.28,

$$\begin{aligned}\mu_{\mathcal{F}'} &\geq 1 + \prod_{x\in N\setminus U}\nu(\mathcal{E}(\{a\},x)) + \prod_{x\in N\setminus U}\nu(\mathcal{E}(\{b\},x)) \\ &\geq 1 + \prod_{x\in N_b}(1-(\tfrac{1}{2})^{n(x)}) + \prod_{x\in N_a}(1-(\tfrac{1}{2})^{n(x)}) \\ &\geq 1 + \left(1-(\tfrac{1}{2})^{n_a}\right)^{n_b} + \left(1-(\tfrac{1}{2})^{n_b}\right)^{n_a}.\end{aligned}$$

Let $c_1 = \left(1-(\tfrac{1}{2})^{n_a}\right)^{n_b}$ and $c_2 = \left(1-(\tfrac{1}{2})^{n_b}\right)^{n_a}$. It remains to show that $1 + c_1 + c_2 \geq 2$. Without loss of generality, assume that $n_a \leq n_b$. If $n_a = 1$, then
$$1 + c_1 + c_2 = 1 + (\tfrac{1}{2})^{n_b} + 1 - (\tfrac{1}{2})^{n_b} = 2.$$

Suppose that $n_a \geq 2$. Observe that $(\frac{3}{2})^m \geq m$ for $m \geq 1$. Using the fact that $(1-a)^b \geq 1 - ba$ for $0 \leq a \leq 1$ and $1 \leq b$ (Lemma 3.30), we obtain

$$\begin{aligned} 1 + c_1 + c_2 &\geq 2 + (\tfrac{3}{4})^{n_b} - n_a (\tfrac{1}{2})^{n_b} \\ &= 2 + (\tfrac{1}{2})^{n_b} \left( (\tfrac{3}{2})^{n_b} - n_a \right) \\ &\geq 2 + (\tfrac{1}{2})^{n_b} \left( (\tfrac{3}{2})^{n_a} - n_a \right) \\ &\geq 2, \end{aligned}$$

as required. ∎

**Lemma 3.30** *If $0 \leq a \leq 1$ and $1 \leq b$, then $(1-a)^b \geq 1 - ba$.*

**Proof.** The result is true for $a = 0$. Differentiating both expressions relative to $a$ yields
$$\frac{d}{da}(1-a)^b = -b(1-a)^{b-1} \geq -b = \frac{d}{da}(1-ba).$$
The result follows. ∎

If $J'(\mathcal{F}) = \emptyset$, then $\mathcal{F}$ is generated by one-element sets, so that $\mathcal{F}$ is a Boolean lattice. Thus by Observation 3.20, Theorem 3.1 has been proved. In fact we have the following stronger result:

**Theorem 3.31** *Let $\mathcal{F}$ be a union-closed family of sets and $U \in J(\mathcal{F})$. Let $\mathcal{G} = \{V \in J(\mathcal{F}) \mid V \cap U \neq \emptyset\}$. If*

(i) $|U| = 2$,

(ii) $\mathcal{G}$ is a graph,

(iii) there is a simple graph $\mathcal{G}' \subseteq J(\mathcal{F})$ such that for every $V \in \mathcal{G}$ with $|V| = 2$, $V \in \mathcal{G}'$ and $d_{\mathcal{G}'}(V) \geq d_{\mathcal{G}'}(U)$,

*then $(\mathcal{F}, U)$ has the density property.*

**Proof.** Suppose that $\mathcal{N}$ is a union-closed family such that $\mathcal{F} \supseteq \mathcal{N} \supseteq \mathcal{N}_\mathcal{F}(U)$. Then $\mathcal{F}$ is an extension of $(\mathcal{N}, U)$ (see the proof of Lemma 3.11). Lemma 3.10 shows that $E_{\mathcal{F},U}(X) \supseteq E_{\mathcal{N},U}(X)$ for every $X$ with $X \cap U = \emptyset$. Therefore

$$\mu_{\mathcal{F},\mathcal{F}'}(U) \geq \mu_{\mathcal{N},\mathcal{F}'}(U)$$

for every extension $\mathcal{F}'$ of $(\mathcal{F}, U)$. Let $\mathcal{N}$ be the union-closed family generated by $\mathcal{G} \cup \mathcal{G}' \cup \{\emptyset\}$. Then $\mathcal{F} \supseteq \mathcal{N} \supseteq \mathcal{N}_\mathcal{F}(U)$. Theorems 3.19 and 3.29 show that $\mu_{\mathcal{N},\mathcal{F}'}(U) \geq 2$ for every extension $\mathcal{F}'$ of $(\mathcal{F}, U)$. It follows that $(\mathcal{F}, U)$ has the density property. ∎